\documentclass{amsart}
\usepackage{amssymb}
\usepackage{amsmath}
\usepackage{amsfonts}

\newtheorem{theorem}{Theorem}
\theoremstyle{plain}

\newtheorem{corollary}{Corollary}

\newtheorem{definition}{Definition}
\newtheorem{example}{Example}

\newtheorem{lemma}{Lemma}
\newtheorem{notation}{Notation}

\numberwithin{equation}{section}
\input{tcilatex}

\begin{document}
\title{Weak conditions for random fixed point and approximation results }
\author{Monica Patriche}
\maketitle

\bigskip\ \ \ \ \ \ \ \ \ \ \ \ \ \ \ \ \ \ \ \ \ \ \ \ \ \ \ \ \ \ \ \ \ \ 
\textit{University of Bucharest}

\begin{center}
\textit{Faculty of Mathematics and Computer Science}

\textit{\ 14 Academiei Street}

\textit{010014 Bucharest, Romania}{\footnotesize \\[0pt]
E-mail: \textit{monica.patriche@yahoo.com}\\[0pt]
\bigskip }
\end{center}

\textbf{Abstract.}{\small \ }In this paper, we study the existence of the
random approximations and fixed points for random almost lower
semicontinuous operators defined on finite dimensional Banach spaces, which
in addition, are condensing or 1-set-contractive. Our results either extend
or improve corresponding ones present in literature.



\section{\textbf{INTRODUCTION}}

Fixed point theorems are crucial in applied mathematics. Various classes of
random equations are based on the random operators theory, which has in its
central part the study of the random fixed points. This study was initiated
by Prague school of probabilists in 1950s and the publication of
Bharucha-Reid's survey paper \cite{b} has been followed by an increasing
interest in this topic. Nowadays, the researchers in this area have devoted
a large amount of activity to prove the existence of the random fixed points
for correspondences. The reader is referred, for instance, to the new
results established in [8], [9], [13], [30], [33]. The case of the
condensing or 1-set-contractive maps has been considered, for example, in
[15], [16], [25], [27], [29], [31], [34].

It is of great interest to obtain random fixed point theorems for
correspondences which verify assumptions weaker than lower semicontinuity.
In this paper, we consider almost lower semicontinuous correspondences, as
well as lower semicontinuous ones, which are, in addition, condensing or
1-set-contractive and we prove approximations and fixed point results for
them. The obtained approximation theorems are stochastic generalizations of
a theorem of Fan [\cite{fan}, Theorem 2], which is stated as follows.

Let $K$ be a nonempty compact convex set in a normed linear space $X$. For
any continuous map $f$ from $K$ into $X$, there exists a point $u$ in $K$
such that, $\left\Vert u-f(x)\right\Vert =d(f(x),K).$

Various versions of this theorem (for single-valued maps or for
correspondences) have been established in the last decades. The reader is
referred, for example, to [12]-[18], [21]-[23], [28], [32], [34].

We mention that, until now, there is no study concerning the existence of
the random approximations in the case of assumptions weaker than lower
semicontinuity. Our research is meant to fill this lack in the theory of
random operators. As applications of our approximation theorems, some
stochastic fixed point theorems are derived.

New research on the existence of fixed points for maps defined on metric
spaces has beed done, for instance, in [3],[4],[26].

The rest of the paper is organized as follows. In the following section,
some notational and terminological conventions are given. We also present,
for the reader's convenience, some information on the continuity and
measurability of the operators. The random approximation and fixed point
theorems for almost lower semicontinuous operators are stated in Section 3.
Section 4 presents the conclusions of our research.

\section{\textbf{NOTATION AND DEFINITION}}

Throughout this paper, we shall use the following notation:

$2^{D}$ denotes the set of all non-empty subsets of the set $D$. If $%
D\subset Y$, where $Y$ is a topological space, cl$D$ denotes the closure of $%
D$. We also denote $C(Y)$ the family of all non-empty and closed subsets of $%
Y.$ A \textit{paracompact} space is a Hausdorff topological space in which
every open cover admits an open locally finite refinement. Metrizable and
compact topological spaces are paracompact.\textit{\ }

Let $X$, $Y$ be topological spaces and $T:X\rightarrow 2^{Y}$ be a
correspondence. $T$ is said to be \textit{lower semicontinuous} if, for each 
$x\in X$ and each open set $V$ in $Y$ with $T(x)\cap V\neq \emptyset $,
there exists an open neighborhood $U$ of $x$ in $X$ such that $T(y)\cap
V\neq \emptyset $ for each $y\in U$.

Let $(X,d)$ be a metric space. We will use the following notations. We
denote by $B(x,r)=\{y\in X:d(y,x)<r\}.$ If $B_{0}$ is a subset of $X,$ then,
we will denote $B(B_{0},r)=\{y\in X:d(y,B_{0})<r\},$ where $%
d(y,B_{0})=\inf_{x\in B_{0}}d(y,x).$

Let $C$ be a \ non-empty subset of $X$ and $T:C\rightarrow 2^{X}$ be a
correspondence.

We say that $T$ is \textit{hemicompact} if each sequence $\{x_{n}\}$ in $C$
has a convergent subsequence, whenever $d(x_{n},T(x_{n}))\rightarrow 0$ as $%
n\rightarrow \infty .$

The correspondence $T:C\rightarrow 2^{X}$ is said to be condensing (see \cite%
{tar}), if for each subset $A$ of $X$ such that $\gamma (A)>0,$ one has $%
\gamma (T(A))<\gamma (A),$ where $T(A)=\cup _{x\in A}T(x)$ and $\gamma $ is
the Kuratowski measure of noncompactness, i.e., for each bounded subset $A$
of $X,$

$\gamma (A)=$inf$\{e>0:A$ is covered by a finite number of sets of diameter $%
\leq e\}.$

If $A$ is not a bound subset of $X,$ we assign $\gamma (A)=\infty .$

If $X$ is a Banach space, the following conditions hold for any $A,B\subset
X $:

(1) $\gamma (A)=0$ if and only if $A$ is precompact;

(2) $\gamma ($co$A)=\gamma (A)$, where co$A$ denotes the closed convex hull
of $A$;

(3) $\gamma (A\cup B)=$max$\{\gamma (A),\gamma (B)\}$.

From (3), we conclude that if $A\subset B,$ then $\gamma (A)\leq \gamma (B)$.

The correspondence $T:C\rightarrow 2^{X}$ is said to be k-set-contractive (%
\cite{lin1}), if for each subset $A$ of $X$ such that $\gamma (A)>0,$ one
has $\gamma (T(A))\leq k\gamma (A).$\smallskip

Let now $(\Omega ,\tciFourier $, $\mu )$ be a complete, finite measure
space, and $Y$ be a topological space. The correspondence $T:\Omega
\rightarrow 2^{Y}$ is said to be \textit{lower measurable} if, for every
open subset $V$ of $Y$, the set $T^{-1}(V)=\{\omega \in \Omega $ $:$ $%
T(\omega )\cap V\neq \emptyset $\} is an element of $\tciFourier $. This
notion of measurability is also called in literature \textit{weak
measurability} or just \textit{measurability}, in comparison with strong
measurability: the correspondence $T:\Omega \rightarrow 2^{Y}$ is said to be 
\textit{strong measurable} if, for every closed subset $V$ of $Y$, the set $%
\{\omega \in \Omega $ $:$ $T(\omega )\cap V\neq \emptyset $\} is an element
of $\tciFourier $. In the case when $Y$ is separable, the strong
measurability coincides with the lower measurability.

Recall (see Debreu \cite{deb2}, p. 359) that if $T:\Omega \rightarrow 2^{Y}$
has a measurable graph, then $T$ is lower measurable. Furthermore, if $%
T(\cdot )$ is closed valued and lower measurable, then $T:\Omega \rightarrow
2^{Y}$ has a measurable graph.

A mapping $T:\Omega \times X\rightarrow Y$ is called a \textit{random
operator} if, for each $x\in X$, the mapping $T(\cdot ,x):\Omega \rightarrow
Y$ is measurable. Similarly, a correspondence $T:\Omega \times X\rightarrow
2^{Y}$ is also called a random operator if, for each $x\in X$, $T(\cdot
,x):\Omega \rightarrow 2^{Y}$ is measurable.

We say that the operator $T:\Omega \times X\rightarrow 2^{Y}$ is condensing
if, for each $\omega \in \Omega ,$ the correspondence $T(\omega ,\cdot
):X\rightarrow 2^{Y}$ is condensing. We say that the operator $T:\Omega
\times X\rightarrow 2^{Y}$ is k-set-contractive if, for each $\omega \in
\Omega ,$ the correspondence $T(\omega ,\cdot ):X\rightarrow 2^{Y}$ is
k-set-contractive.

A measurable mapping $\xi :\Omega \rightarrow Y$ is called a \textit{%
measurable selection of the operator} $T:\Omega \rightarrow 2^{Y}$ if $\xi
(\omega )\in T(\omega )$ for each $\omega \in \Omega $. A measurable mapping 
$\xi :\Omega \rightarrow Y$ is called a \textit{random fixed point} of the
random operator $T:\Omega \times X\rightarrow Y$ (or $T:\Omega \times
X\rightarrow 2^{Y})$ if for every $\omega \in \Omega $%
\c{}
$\xi (\omega )=T(\omega ,\xi (\omega ))$ (or $\xi (\omega )\in T(\omega ,\xi
(\omega ))$).

We will need the following measurable selection theorem in order to prove
our results. \smallskip

\textbf{Proposition 2.1} (Kuratowski-Ryll-Nardzewski Selection Theorem \cite%
{k}). A weakly measurable correspondence with non-empty closed values from a
measurable space into a Polish space admits a measurable selector.

\section{RANDOM FIXED POINT AND APPROXIMATION THEOREMS FOR RANDOM ALMOST
LOWER SEMICONTINUOUS OPERATORS}

This section is mainly dedicated to establishing the random approximation
and fixed point results concerning the random almost lower semicontinuous
condensing (or\textit{\ }1-set-contractive\textit{) }operators.

Firstly we recall the following statement, which will be useful to prove the
first result of this section.\smallskip

\begin{lemma}
\textit{(Theorem 3.4.in \cite{fie1}) Let }$C$\textit{\ be a closed,
separable subset of a complete metric space }$X$\textit{, and let }$T:%
\Omega
\times C\rightarrow C(X)$\textit{\ be a continuous hemicompact random
operator. If, for each }$\omega \in 
\Omega
,$\textit{\ the set }$F(\omega ):=\{x\in C:x\in T(\omega ,x)\}\neq \emptyset 
$\textit{, then, }$T$\textit{\ has a random fixed point.\smallskip }
\end{lemma}

In order to extend the fixed point theorems to condensing operators, we will
use Lemma 2 (see \cite{m}).

\begin{lemma}
\cite{m} Let $X$ denote a nonempty, closed and convex subset of a Hausdorff
locally convex topological vector space $E$. If $T:X\rightarrow 2^{X}$ is
condensing, then there exists a nonempty, compact and convex subset $K$ of $%
X $ such that $T(x)\subset K$ for each $x\in K.\smallskip $
\end{lemma}

Now, we are presenting the almost lower semicontinuous correspondences.

Let $X$\ be a topological space and $Y$\ be a normed linear space.\textit{\ }%
The correspondence $T:X\rightarrow 2^{Y}$ is said to be \textit{almost lower
semicontinuous (a.l.s.c.)} \textit{at} $x\in X$ (see \cite{deu}), if, for
any $\varepsilon >0,$ there exists a neighborhood $U(x)$ of $x$ such that $%
\tbigcap\limits_{z\in U(x)}B(T(z);\varepsilon )\neq \emptyset .$

T is \textit{almost lower semicontinuous} if it is a.l.s.c. at each $x\in X$.

If $\Omega $ is a non-empty set, we say that the operator $T:\Omega \times
X\rightarrow 2^{Y}$ is almost lower semicontinuous if, for each $\omega \in
\Omega ,$ $T(\omega ,\cdot )$ is \ almost lower semicontinuous.

In 1983, Deutsch and Kenderov \cite{deu} presented a remarkable
characterization of a.l.s.c. correspondences as follows.\smallskip

\begin{lemma}
\textit{(see \cite{deu}) Let }$X$\textit{\ be a paracompact topological
space, }$Y$\textit{\ be a normed vector space and }$T:X\rightarrow 2^{Y}$%
\textit{\ be a correspondence having convex values. Then, }$\mathit{T}$%
\textit{\ is a.l.s.c. if and only if, for each }$\varepsilon >0$\textit{, }$%
T $\textit{\ admits a continuous }$\varepsilon -$\textit{approximate
selection f; that is, }$f:X\rightarrow Y$\textit{\ is a continuous
single-valued function such that }$f(x)\in B(T(x);\varepsilon )$\textit{\
for each }$x\in X.\smallskip $
\end{lemma}

The next theorem states the existence of random fixed points for the random
almost lower semicontinuous condensing operators defined on Banach
spaces.\smallskip

\begin{theorem}
\textit{Let }$(\Omega ,\mathcal{F})$\textit{\ be a measurable space, }$C$%
\textit{\ be a closed convex separable subset of a finite dimensional Banach
space }$X$\textit{\ and let }$T:\Omega \times C\rightarrow 2^{C}$\textit{\
be a random operator. Suppose that, for each }$\omega \in \Omega ,$\textit{\ 
}$T(\omega ,\cdot )$\textit{\ is almost lower semicontinuous, condensing,
with non-empty, convex and closed values. In addition, assume that there is
a }$n_{0}\in \mathbb{N}^{\ast }$ with the property that $B(T(\omega ,C),%
\frac{1}{n_{0}})\subseteq C$ \textit{and }$(T(\omega ,\cdot
))^{-1}:C\rightarrow 2^{C}$\textit{\ is closed valued.}
\end{theorem}

\textit{Then, }$T$\textit{\ has a random fixed point.}

\begin{proof}
If for each $\omega \in \Omega ,$ $T(\omega ,\cdot ):C\rightarrow 2^{C}$ is
condensing, then $T_{n_{0}}(\omega ,\cdot ):C\rightarrow 2^{C}$ defined by $%
T_{n_{0}}(\omega ,x)=B(T(\omega ,x);1/n_{0})$ if $(\omega ,x)\in \Omega
\times C$ is also condensing. Indeed, for each bounded subset $A$ of $C,$ $%
\gamma (T_{n_{0}}(\omega ,A))=\gamma (T(\omega ,A)+B(0,\frac{1}{n_{0}}))\leq
\gamma (T(\omega ,A))+\gamma (B(0,\frac{1}{n_{0}}))=\gamma (T(\omega
,A))<\gamma (A).$ (The last inequality is true when $T$ is condensing). We
note that $\gamma (A)=0$ if and only if $A$ is precompact and $B(0,\frac{1}{%
n_{0}})$ is precompact in $X,$ and therefore, \ $\gamma (B(0,\frac{1}{n_{0}}%
))=0.$

According to Lemma 2, if $T_{n_{0}}(\omega ,\cdot ):C\rightarrow 2^{C}$ is
condensing, then there exists a nonempty compact convex subset $K$ of $C$
such that $T_{n_{0}}(\omega ,x)\subset K$ for each $x\in K.$

Firstly, let us define $T_{n}(\omega ,\cdot ):K\rightarrow 2^{K}$ by $%
T_{n}(\omega ,x)=B(T(\omega ,x);1/(n+n_{0}-1))$ if $(\omega ,x)\in \Omega
\times K$ and $n\in \mathbb{N}^{\ast }.$ Since for each $\omega \in \Omega ,$
$T(\omega ,\cdot )$ is almost lower semicontinuous, according to Lemma 3,
for each $n\in N,$\ there exists a continuous function $f_{n}(\omega ,\cdot
):K\rightarrow K$\ such that $f_{n}(\omega ,x)\in T_{n}(\omega ,x)$ for each%
\textit{\ }$x\in K.$ Brouwer-Schauder fixed point theorem enssures that, for
each $n\in N,$ there exists $x_{n}\in K$ such that $x_{n}=f_{n}(\omega
,x_{n})$ and then, $x_{n}\in T_{n}(\omega ,x_{n}).$

$K$ is compact, so $f_{n}$ is hemicompact for each $n\in \mathbb{N}$.
According to Lemma 1, for each $n\in \mathbb{N},$ $f_{n}$ has a random fixed
point and then, $T_{n}$ has a random fixed point $\xi _{n},$ that is $\xi
_{n}:\Omega \rightarrow K$ is measurable and $\xi _{n}(\omega )\in
T_{n}(\omega ,\xi _{n}(\omega ))$ for $n\in N$.

Let $\omega \in \Omega $ be fixed. Then, $d(\xi _{n}(\omega ),T(\omega ,\xi
_{n}(\omega ))\rightarrow 0$ when $n\rightarrow \infty $ and since $K$ is
compact, $\{\xi _{n}(\omega )\}$ has a convergent subsequence $\{\xi
_{n_{k}}(\omega )\}.$ Let $\xi _{0}(\omega )=\lim_{n_{k}\rightarrow \infty
}\xi _{n_{k}}(\omega ).$ It follows that $\xi _{0}:\Omega \rightarrow K$ is
measurable and for each $\omega \in \Omega ,$ $d(\xi _{0}(\omega ),T(\omega
,\xi _{n_{k}}(\omega ))\rightarrow 0$ when $n_{k}\rightarrow \infty .$

Let us assume that there is a $\omega \in \Omega $ such that $\xi
_{0}(\omega )\notin T(\omega ,\xi _{0}(\omega )).$ Since $\{\xi _{0}(\omega
)\}\cap (T(\omega ,\cdot ))^{-1}(\xi _{0}(\omega ))=\emptyset $ and $X$ is a
regular space, there exists $r_{1}>0$ such that $B(\xi _{0}(\omega
),r_{1})\cap (T(\omega ,\cdot ))^{-1}(\xi _{0}(\omega ))=\emptyset $.
Consequently, for each $z\in B(\xi _{0}(\omega ),r_{1}),$ we have that $%
z\notin (T(\omega ,\cdot ))^{-1}(\xi _{0}(\omega )),$ which is equivalent
with $\xi _{0}(\omega )\notin T(\omega ,z)$ or $\{\xi _{0}(\omega )\}\cap
T(\omega ,z)=\emptyset $. The closedness of each $T(\omega ,z)$ and the
regularity of $X$ imply the existence of a real number $r_{2}>0$ such that $%
B(\xi _{0}(\omega ),r_{2})\cap T(\omega ,z)=\emptyset $ for each $z\in B(\xi
_{0}(\omega ),r_{1}),$ which implies $\xi _{0}(\omega )\notin B(T(\omega
,z);r_{2})$ for each $z\in B(\xi _{0}(\omega ),r_{1}).$ Let $r=\min
\{r_{1},r_{2}\}.$ Hence, $\xi _{0}(\omega )\notin B(T(\omega ,z);r)$ for
each $z\in B(\xi _{0}(\omega ),r),$ and then, there exists $N^{\ast }\in 
\mathbb{N}$ such that for each $n_{k}>N^{\ast },$ $\xi _{0}(\omega )\notin
B(T(\omega ,\xi _{n_{k}}(\omega ));r)$ which contradicts $d(\xi _{0}(\omega
),T(\omega ,\xi _{n_{k}}(\omega ))\rightarrow 0$ as $n\rightarrow \infty $.
It follows that our assumption is false.

Hence, for each $\omega \in \Omega ,$ $\xi _{0}(\omega )\in T(\omega ,\xi
_{0}(\omega ))$, where $\xi _{0}:\Omega \rightarrow K$ is measurable. We
conclude that $T$ has a random fixed point.
\end{proof}

\begin{example}
Let $T_{1}:[0,\infty )\rightarrow 2^{[0,\infty )}$ be defined by
\end{example}

$T_{1}(x)=\left\{ 
\begin{array}{c}
\{0.00005\},\text{ \ \ \ \ \ \ if \ \ \ \ \ \ }x\in \lbrack 0,\frac{1}{100}];
\\ 
\{\frac{1}{2}x^{2}\},\text{ \ \ if \ \ }x\in (\frac{1}{100},\frac{15}{32}%
)\cup (\frac{15}{32},1]; \\ 
\lbrack \frac{1}{10},\frac{1}{2}],\text{ \ \ \ \ \ \ \ \ \ \ \ if \ \ \ \ \
\ \ \ \ \ \ \ }x=\frac{15}{32}; \\ 
\{\frac{1}{2}\},\text{ \ \ \ \ \ \ \ \ \ \ if \ \ \ \ \ \ \ \ \ \ \ \ \ \ \
\ \ }x>1.%
\end{array}%
\right. $

\textit{We note that} $T_{1}$\textit{\ is almost lower semicontinuous,
condensing, with non-empty, convex and closed values.}

\textit{Let }$\Omega =[0,\infty ),$\textit{\ }$F$ \textit{be} \textit{the }$%
\sigma -$\textit{algebra of the borelian sets of }$[0,\infty )$\textit{\ and
let }$T:\Omega \times \lbrack 0,\infty )\rightarrow 2^{[0,\infty )}$\textit{%
\ be the random operator defined by}

$T(\omega ,x)=\left\{ 
\begin{array}{c}
T_{1}(x)\text{ \ \ if \ \ \ \ \ }x=\omega ; \\ 
\lbrack 0,00005,\frac{1}{2}]\text{ if }x\neq \omega%
\end{array}%
\right. $ for each $(\omega ,x)\in \Omega \times \lbrack 0,\infty ).$

\textit{For each }$\omega \in \Omega ,$\textit{\ }$T(\omega ,\cdot )$\textit{%
\ is almost lower semicontinuous, condensing, with non-empty, convex and
closed values. There exists }$n_{0}=20000\in \mathbb{N}^{\ast }$ \textit{%
with the property that} $B(T(\omega ,[0,\infty )),\frac{1}{n_{0}})\subseteq
\lbrack 0,\infty ).$

\textit{We will prove that, for each }$\omega \in \Omega ,$ $(T(\omega
,\cdot ))^{-1}:[0,\infty )\rightarrow 2^{[0,\infty )}$\textit{\ is closed
valued.}

\textit{If }$\omega \in \lbrack 0,\frac{1}{100}],$

$T^{-1}(\omega ,y)=\left\{ 
\begin{array}{c}
\phi ,\text{ \ \ \ \ \ \ \ \ \ \ \ \ \ \ \ if \ \ \ \ \ \ \ \ \ \ \ \ \ \ \
\ \ \ \ \ \ \ }y\in \lbrack 0,0.00005); \\ 
\lbrack 0,\frac{1}{100}],\text{ \ \ \ \ \ \ if \ \ \ \ \ \ \ \ \ \ \ \ \ \ \
\ \ \ \ \ \ \ \ \ \ \ \ \ \ \ }y=0.00005; \\ 
\{\omega ,\sqrt{2y}\},\text{ \ \ \ \ \ \ \ \ if \ \ \ \ \ \ \ \ \ \ \ \ \ \
\ \ \ }y\in (0.00005,0.10); \\ 
\{\omega ,\sqrt{2y},\frac{15}{32}\},\text{ }\ \ \ \ \ \ \text{if }y\in
\lbrack 0.10,0.109)\cup (0.109,0.5); \\ 
\{\omega ,\frac{15}{32}\},\text{ \ \ \ if \ \ \ \ \ \ \ \ \ \ \ \ \ \ \ \ \
\ \ \ \ \ \ \ \ \ \ \ \ \ \ \ \ \ \ \ }y=0.109; \\ 
\{\omega ,\frac{15}{32}\}\cup \lbrack 1,\infty ),\text{\ \ \ \ if \ \ \ \ \
\ \ \ \ \ \ \ \ \ \ \ \ \ \ \ \ \ \ \ \ \ \ }y=0.5; \\ 
\phi ,\text{ \ \ \ \ \ \ \ \ \ if \ \ \ \ \ \ \ \ \ \ \ \ \ \ \ \ \ \ \ \ \
\ \ \ \ \ \ \ \ \ \ \ }y\in (0.5,\infty ).%
\end{array}%
\right. $

\textit{If }$\omega \in (\frac{1}{100},\frac{15}{32}),$

$T^{-1}(\omega ,y)=\left\{ 
\begin{array}{c}
\phi ,\text{ \ \ \ \ \ \ \ \ if \ \ \ \ \ \ \ \ \ \ \ \ \ \ \ \ \ \ \ \ \ \
\ \ }y\in \lbrack 0,0.00005); \\ 
\lbrack 0,\frac{1}{100}]\cup \{\omega \},\text{ \ \ \ \ \ \ \ \ \ if \ \ \ \
\ \ \ \ \ \ \ \ }y=0.00005; \\ 
\{\omega ,\sqrt{2y}\},\text{ }\ \ \ \ \ \ \ \ \text{if \ \ \ \ \ \ \ \ \ \ }%
y\in (0.00005,0.10); \\ 
\{\omega ,\sqrt{2y},\frac{15}{32}\},\text{ if }y\in \lbrack 0.10,0.109)\cup
(0.109,0.5); \\ 
\{\omega ,\frac{15}{32}\},\text{ \ \ \ \ \ \ \ \ \ \ \ if \ \ \ \ \ \ \ \ \
\ \ \ \ \ \ \ \ \ \ \ \ }y=0.109; \\ 
\{\omega ,\frac{15}{32}\}\cup \lbrack 1,\infty ),\text{ \ \ \ \ \ \ \ \ if \
\ \ \ \ \ \ \ \ \ \ \ \ \ \ }y=0.5; \\ 
\phi ,\text{ \ \ \ \ \ \ \ if \ \ \ \ \ \ \ \ \ \ \ \ \ \ \ \ \ \ \ \ \ \ \
\ \ \ \ }y\in (0.5,\infty ).%
\end{array}%
\right. $

\textit{If }$\omega =\frac{15}{32},$

$T^{-1}(\omega ,y)=\left\{ 
\begin{array}{c}
\phi ,\text{ \ \ \ \ \ \ \ \ \ \ \ \ \ \ \ \ \ \ if \ \ \ \ \ \ \ \ \ }y\in
\lbrack 0,0.00005); \\ 
\lbrack 0,\frac{1}{100}]\cup \{\frac{15}{32}\},\text{ \ \ \ \ \ \ \ \ if \ \
\ \ \ \ \ \ }y=0.00005; \\ 
\{\frac{15}{32},\sqrt{2y}\},\text{ \ \ if \ \ }y\in (0.00005,0.5)\backslash
\{0.109\}; \\ 
\{\frac{15}{32}\},\text{ \ \ \ \ \ \ \ \ \ \ \ \ if \ \ \ \ \ \ \ \ \ \ \ \
\ \ \ \ \ \ \ \ \ }y=0.109; \\ 
\{\frac{15}{32}\}\cup \lbrack 1,\infty ),\text{ \ \ \ \ \ \ \ \ \ \ \ \ \ \
\ \ if \ \ \ \ \ \ \ \ \ }y=0.5; \\ 
\phi ,\text{ \ \ \ \ \ \ \ \ \ \ \ \ \ \ \ \ \ if \ \ \ \ \ \ \ \ \ \ \ \ \ }%
y\in \lbrack 0.05,\infty ).%
\end{array}%
\right. $

\textit{If }$\omega \in (\frac{15}{32},1),$

$T^{-1}(\omega ,y)=\left\{ 
\begin{array}{c}
\phi ,\text{ \ \ \ \ \ \ \ \ \ if \ \ \ \ \ \ \ \ \ \ \ \ \ \ \ \ \ \ \ \ \
\ }y\in \lbrack 0,0.00005); \\ 
\lbrack 0,\frac{1}{100}]\cup \{\omega \},\text{ \ \ \ \ \ \ \ \ \ \ \ \ \ \
\ \ \ if \ \ \ \ }y=0.00005; \\ 
\{\omega ,\sqrt{2y}\},\text{ \ \ \ \ \ \ \ \ \ \ \ \ if \ \ \ \ \ \ \ }y\in
(0.00005,0.10); \\ 
\{\omega ,\sqrt{2y},\frac{15}{32}\},\text{ if }y\in \lbrack 0.10,0.109)\cup
(0.109,0.5); \\ 
\{\omega ,\frac{15}{32}\},\text{ \ \ \ \ \ \ \ \ \ \ \ \ \ \ \ if \ \ \ \ \
\ \ \ \ \ \ \ \ \ \ \ \ \ }y=0.109; \\ 
\{\omega ,\frac{15}{32}\}\cup \lbrack 1,\infty ),\text{ \ \ \ \ \ \ \ \ \ \
\ \ \ if \ \ \ \ \ \ \ \ \ \ \ \ }y=0.5; \\ 
\phi ,\text{ \ \ \ \ \ \ \ \ if \ \ \ \ \ \ \ \ \ \ \ \ \ \ \ \ \ \ \ \ \ \
\ \ \ }y\in (0.05,\infty ).%
\end{array}%
\right. $

\textit{If }$\omega \in \lbrack 1,\infty ),$

$T^{-1}(\omega ,y)=\left\{ 
\begin{array}{c}
\phi ,\text{ \ \ \ \ if \ \ \ \ \ \ \ \ \ \ \ \ \ \ \ \ \ \ \ \ \ \ \ \ \ \
\ }y\in \lbrack 0,0.00005); \\ 
\lbrack 0,\frac{1}{100}]\cup \{\omega \},\text{ \ \ \ \ \ \ \ \ \ \ if \ \ \
\ \ \ \ \ \ \ \ }y=0.00005; \\ 
\{\omega ,\sqrt{2y}\},\text{ \ \ \ \ \ \ \ if \ \ \ \ \ \ \ \ \ \ \ }y\in
(0.00005,0.10); \\ 
\{\omega ,\sqrt{2y},\frac{15}{32}\},\text{ if }y\in \lbrack 0.10,0.109)\cup
(0.109,0.5); \\ 
\{\omega ,\frac{15}{32}\},\text{ \ \ \ \ \ \ \ \ \ \ \ \ \ \ \ \ if \ \ \ \
\ \ \ \ \ \ \ \ \ \ \ \ \ }y=0.109; \\ 
\{\frac{15}{32}\}\cup \lbrack 1,\infty ),\text{ \ \ \ \ \ \ \ \ \ if \ \ \ \
\ \ \ \ \ \ \ \ \ \ \ \ \ \ \ }y=0.5; \\ 
\phi ,\text{ \ \ if \ \ \ \ \ \ \ \ \ \ \ \ \ \ \ \ \ \ \ \ \ \ \ \ \ \ \ \
\ \ \ \ }y\in (0.05,\infty ).%
\end{array}%
\right. $

\textit{The relations written above prove that} \textit{for each }$\omega
\in \Omega ,$ $(T(\omega ,\cdot ))^{-1}:[0,\infty )\rightarrow 2^{[0,\infty
)}$\textit{\ is closed valued.}

\textit{All the conditions of Theorem 1 are fulfilled and thus, there exists 
}$\xi :\Omega \rightarrow \lbrack 0,\infty )$\textit{\ a measurable function
such that }$\xi (\omega )\in T(\omega ,\xi (\omega )).$

Let $\xi :\Omega \rightarrow \lbrack 0,\infty )$ be defined by $\xi (\omega
)=0.00005$ for each $\omega \in \Omega .$

$T(\omega ,\xi (\omega ))=\left\{ 
\begin{array}{c}
T_{1}(\xi (\omega )),\text{ \ if \ \ \ }\xi (\omega )=\omega ; \\ 
\lbrack 0.00005,\frac{1}{2}]\text{ if }\xi (\omega )\neq \omega .%
\end{array}%
\right. =$

$=\left\{ 
\begin{array}{c}
T_{1}(0.00005),\text{ \ \ \ \ \ \ \ \ \ \ \ \ \ if \ \ \ \ \ \ \ \ \ \ \ }%
\omega =0.00005; \\ 
\lbrack 0.00005,\frac{1}{2}]\text{ if }\omega \in \lbrack 0,0.00005)\cup
(0.00005,\infty ).%
\end{array}%
\right. $

$=\left\{ 
\begin{array}{c}
\{0.00005\},\text{ \ \ \ \ \ \ \ \ \ \ \ \ \ if \ \ \ \ \ \ \ \ \ \ \ \ }%
\omega =0.00005; \\ 
\lbrack 0.00005,\frac{1}{2}]\text{ if }\omega \in \lbrack 0,0.00005)\cup
(0.00005,\infty ).%
\end{array}%
\right. $

Therefore, $\xi (\omega )=0.00005\in T(\omega ,\xi (\omega ))$ for each $%
\omega \in \Omega .$

\begin{corollary}
\textit{Let }$(\Omega ,\mathcal{F})$\textit{\ be a measurable space, }$C$%
\textit{\ be a closed convex separable subset of a finite dimensional Banach
space }$X$\textit{\ and let }$T:\Omega \times C\rightarrow 2^{C}$\textit{\
be a random operator. Suppose that, for each }$\omega \in \Omega ,$\textit{\ 
}$T(\omega ,\cdot )$\textit{\ is lower semicontinuous, condensing, with
non-empty, convex and closed values. In addition, assume that there exists }$%
n_{0}\in \mathbb{N}^{\ast }$ with the property that $B(T(\omega ,C),\frac{1}{%
n_{0}})\subseteq C$ \textit{and }$(T(\omega ,\cdot ))^{-1}:C\rightarrow
2^{C} $\textit{\ is closed valued.}
\end{corollary}

\textit{Then, }$T$\textit{\ has a random fixed point.}

\begin{notation}
We denote $B_{R}=\{x\in X:\left\Vert x\right\Vert \leq R\}$, $\partial
B_{1}=\{x\in X:\left\Vert x\right\Vert =1\}$ and $E(X,B_{1})=\{B(x,r)\subset
X:x\in B_{1},r\geq 0\}.$
\end{notation}

\ \ \ \ \ \ Let $C$ be a subset of a Hausdorff topological vector space $X$
and $x\in X.$ Then the inward set $I_{C}(x)$ is defined by

$I_{C}(x)=\{x+r(y-x):y\in C,$ $r\geq 0\}.$

If $C$ is convex and $x\in C,$ then,

$I_{C}(x)=x+\{r(y-x):y\in C,$ $r\geq 1\}.\medskip $

Now, we are proving a random approximation theorem for random almost lower
semicontinuous and condensing operators.

\begin{theorem}
\textit{Let }$(\Omega ,\mathcal{F})$\textit{\ be a measurable space, }$B_{2}$%
\textit{\ be separable in a finite dimensional Banach space }$X$\textit{\
and let }$T:\Omega \times B_{2}\rightarrow E(X,B_{1})$ \textit{be a random
operator. Suppose that, for each }$\omega \in \Omega ,$\textit{\ }$T(\omega
,\cdot )$\textit{\ is almost lower semicontinuous, condensing, non-empty
valued and }$(T(\omega ,\cdot ))^{-1}:X\rightarrow 2^{B_{2}}$\textit{\ is
closed valued.}
\end{theorem}

\textit{Then, there exist a measurable mapping }$\xi :\Omega \rightarrow
B_{1}$\textit{\ and a mapping }$\eta :\Omega \rightarrow X$ \textit{such that%
} \textit{for each} $\omega \in \Omega ,$ \textit{we have}

\ \ \ \ \ \ \ \ \ \ \ \ \ \ \ \ \ \ \ \ \ \ \ \ \ \ \ \ \ \ \ \ \ \ \ \ \ \ $%
\eta (\omega )\in T(\omega ,\xi (\omega ))$

\textit{and}

\ \ \ \ \ \ \ \ \ \ \ \ \ \ \ \ \ \ \ \ \ \ \ \ \ \ \ \ \ \ $\left\Vert \eta
(\omega )-\xi (\omega )\right\Vert =d(\eta (\omega ),B_{1})=d(\eta (\omega ),%
\overline{I_{B_{1}}(\xi (\omega ))}).$

\begin{proof}
Let us define $r:X\rightarrow B_{1}$ by $r(x)=\left\{ 
\begin{array}{c}
x\text{ if }x\in B_{1}; \\ 
\frac{x}{\left\Vert x\right\Vert }\text{ if }x\notin B_{1}.%
\end{array}%
\right. $

Then, $r$ is continuous and $r(A)\subseteq \overline{\text{co}}(A\cup \{0\})$
for each bounded subset $A$ of $X.$ Thus, $\gamma (r(A))\leq \gamma (A)$ and
this means that $r$ is a 1-set-contractive map. In addition, $T(\omega
,\cdot )$ is condensing for each $\omega \in \Omega ,$ and we conclude that $%
G(\omega ,\cdot )=r\circ T(\omega ,\cdot ):B_{2}\rightarrow 2^{B_{1}}$ is
condensing. According to the hypotheses, for each $\omega \in \Omega ,$\ $%
T(\omega ,\cdot )$\ is almost lower semicontinuous with non-empty, convex
and closed values and $(T(\omega ,\cdot ))^{-1}:X\rightarrow 2^{B_{2}}$\ is
closed valued. These imply that for each $\omega \in \Omega ,$\ $G(\omega
,\cdot )$\ is almost lower semicontinuous with non-empty, convex and closed
values and $(G(\omega ,\cdot ))^{-1}:B_{1}\rightarrow 2^{B_{2}}$\ is closed
valued.

For each $n\in \mathbb{N}^{\ast }$ and $\omega \in \Omega ,$ let $%
G_{n}(\omega ,\cdot ):B_{2}\rightarrow 2^{B_{2}}$ be defined by $%
G_{n}(\omega ,x)=B(G(\omega ,x),\frac{1}{n}).$ Each $G_{n}(\omega ,\cdot )$
is also condensing. Indeed, for each bounded subset $A$ of $B_{2},$ $\gamma
(G_{n}(\omega ,A))=\gamma (G(\omega ,A)+B(0,\frac{1}{n}))\leq \gamma
(G(\omega ,A))+\gamma (B(0,\frac{1}{n}))=\gamma (G(\omega ,A))<\gamma (A).$

According to Lemma 2, since $G_{1}(\omega ,\cdot ):B_{2}\rightarrow
2^{B_{2}} $ is condensing, there exists a nonempty, compact and convex
subset $K$ of $B_{2},$ such that $G_{1}(\omega ,x)\subset K$ for each $x\in
K.$

Then, for each $n\in \mathbb{N}^{\ast },$ $G_{n}(\omega ,K)\subset K.$ The
correspondence $G$ is almost lower semicontinuous, and then, according to
Lemma 2, for each $n\in \mathbb{N}^{\ast },$\ there exists a continuous
function $f_{n}(\omega ,\cdot ):K\rightarrow K$\ such that $f_{n}(\omega
,x)\in G_{n}(\omega ,x)$ for each\textit{\ }$x\in K.$ Brouwer-Schauder fixed
point theorem enssures that, for each $n\in N,$ there exists $x_{n}\in K$
such that $x_{n}=f_{n}(\omega ,x_{n})$ and then, $x_{n}\in G_{n}(\omega
,x_{n}).$

$K$ is compact, then $f_{n}$ is hemicompact for each $n\in \mathbb{N}$.
According to Lemma 1, for each $n\in \mathbb{N},$ $f_{n}$ has a random fixed
point and then, $G_{n}$ has a random fixed point $\xi _{n},$ that is, $\xi
_{n}:\Omega \rightarrow K$ is measurable and $\xi _{n}(\omega )\in
G_{n}(\omega ,\xi _{n}(\omega ))$ for $n\in N$.

Let $\omega \in \Omega $ be fixed. Then, $d(\xi _{n}(\omega ),G(\omega ,\xi
_{n}(\omega ))\rightarrow 0$ when $n\rightarrow \infty $ and since $K$ is
compact, $\{\xi _{n}(\omega )\}$ has a convergent subsequence $\{\xi
_{n_{k}}(\omega )\}.$ Let $\xi _{0}(\omega )=\lim_{n_{k}\rightarrow \infty
}\xi _{n_{k}}(\omega ).$ It follows that $\xi _{0}:\Omega \rightarrow K$ is
measurable and for each $\omega \in \Omega ,$ $d(\xi _{0}(\omega ),G(\omega
,\xi _{n_{k}}(\omega ))\rightarrow 0$ when $n_{k}\rightarrow \infty .$

Let us assume that there exists $\omega \in \Omega $ such that $\xi
_{0}(\omega )\notin G(\omega ,\xi _{0}(\omega )).$ Since $\{\xi _{0}(\omega
)\}\cap (G(\omega ,\cdot ))^{-1}(\xi _{0}(\omega ))=\emptyset $ and $X$ is a
regular space, there exists $r_{1}>0$ such that $B(\xi _{0}(\omega
),r_{1})\cap (G(\omega ,\cdot ))^{-1}(\xi _{0}(\omega ))=\emptyset $.
Consequently, for each $z\in B(\xi _{0}(\omega ),r_{1}),$ we have that $%
z\notin (G(\omega ,\cdot ))^{-1}(\xi _{0}(\omega )),$ which is equivalent
with $\xi _{0}(\omega )\notin G(\omega ,z)$ or $\{\xi _{0}(\omega )\}\cap
G(\omega ,z)=\emptyset $. The closedness of each $G(\omega ,z)$ and the
regularity of $X$ imply the existence of a real number $r_{2}>0$ such that $%
B(\xi _{0}(\omega ),r_{2})\cap G(\omega ,z)=\emptyset $ for each $z\in B(\xi
_{0}(\omega ),r_{1}),$ which implies $\xi _{0}(\omega )\notin B(G(\omega
,z);r_{2})$ for each $z\in B(\xi _{0}(\omega ),r_{1}).$ Let $r=\min
\{r_{1},r_{2}\}.$ Hence, $\xi _{0}(\omega )\notin B(G(\omega ,z);r)$ for
each $z\in B(\xi _{0}(\omega ),r),$ and then, there exists $N^{\ast }\in 
\mathbb{N}$ such that for each $n_{k}>N^{\ast },$ $\xi _{0}(\omega )\notin
B(G(\omega ,\xi _{n_{k}}(\omega ));r)$ which contradicts $d(\xi _{0}(\omega
),G(\omega ,\xi _{n_{k}}(\omega ))\rightarrow 0$ as $n\rightarrow \infty $.
It follows that our assumption is false.

Hence, we obtain that there exists a measurable mapping $\xi :\Omega
\rightarrow B_{1}$ such that $\xi (\omega )\in G(\omega ,\xi (\omega ))$ for
each $\omega \in \Omega .$

Let $\eta :\Omega \rightarrow B_{1}$ be such that $\xi (\omega )=r(\eta
(\omega ))$ for each $\omega \in \Omega .$ Then, $\eta (\omega )\in T(\omega
,\xi (\omega ))$ for each $\omega \in \Omega .$

Let $\omega \in \Omega $ be fixed.

Further, we will consider the cases: $\eta (\omega )\in B_{1}$ and $\eta
(\omega )\notin B_{1}.$

If $\eta (\omega )\in B_{1},$ it is obvious that $\xi (\omega )=r(\eta
(\omega ))=\eta (\omega )$ and consequently, $\left\Vert \eta (\omega )-\xi
(\omega )\right\Vert =0=d(\eta (\omega ),B_{1}).$

In case that $\eta (\omega )\notin B_{1},$ then $\xi (\omega )=r(\eta
(\omega ))=\frac{\eta (\omega )}{\left\Vert \eta (\omega )\right\Vert }.$
This implies that for each $x\in B_{1},$ $\left\Vert \eta (\omega
)-x\right\Vert \geq \left\Vert \eta (\omega )\right\Vert -\left\Vert
x\right\Vert \geq \left\Vert \eta (\omega )\right\Vert -1=\frac{\left\Vert
\eta (\omega )\right\Vert -1}{\left\Vert \eta (\omega )\right\Vert }%
\left\Vert \eta (\omega )\right\Vert =\left\Vert \eta (\omega )-\frac{\eta
(\omega )}{\left\Vert \eta (\omega )\right\Vert }\right\Vert =\left\Vert
\eta (\omega )-\xi (\omega )\right\Vert .$

Therefore, $\left\Vert \eta (\omega )-\xi (\omega )\right\Vert =d(\eta
(\omega ),B_{1})$ for each $\omega \in \Omega .$

We will further prove the equality:

$d(\eta (\omega ),B_{1})=d(\eta (\omega ),\overline{I_{B_{1}}(\xi (\omega ))}%
)$ for each $\omega \in \Omega .$

In order to do this, we choose arbitrarily $\omega \in \Omega $ and we
consider $z\in I_{B_{1}}(\xi (\omega ))\backslash B_{1}.$ There exist $y\in
B_{1}$ and $\lambda >1$ such that $z=\xi (\omega )+\lambda (y-\xi (\omega
)). $

By way of contradiction, we suppose $\left\Vert \eta (\omega )-z\right\Vert
<\left\Vert \eta (\omega )-\xi (\omega )\right\Vert .$

Since $\frac{1}{\lambda }z+(1-\frac{1}{\lambda })\xi (\omega )\in B_{1},$ we
obtain

$\left\Vert \eta (\omega )-z\right\Vert =\left\Vert \frac{1}{\lambda }(\eta
(\omega )-z)+(1-\frac{1}{\lambda })(\eta (\omega )-\xi (\omega ))\right\Vert 
$

\ \ \ \ \ \ \ \ \ \ \ \ \ \ \ \ \ \ $\leq \frac{1}{\lambda }\left\Vert (\eta
(\omega )-z)\right\Vert +(1-\frac{1}{\lambda })\left\Vert (\eta (\omega
)-\xi (\omega ))\right\Vert $

$<\left\Vert (\eta (\omega )-\xi (\omega ))\right\Vert ,$ which is a
contradiction.

Hence,

\ \ \ \ \ \ \ \ \ \ \ \ \ \ \ \ \ \ $\left\Vert \eta (\omega )-\xi (\omega
)\right\Vert \leq \left\Vert \eta (\omega )-z\right\Vert $

for each $z\in I_{B_{1}}(\xi (\omega )),$ and thus we proved that

$\left\Vert \eta (\omega )-\xi (\omega )\right\Vert =d(\eta (\omega
),B_{1})=d(\eta (\omega ),\overline{I_{B_{1}}(\xi (\omega ))})$ for each $%
\omega \in \Omega .$
\end{proof}

\bigskip

\begin{example}
Let $X=\mathbb{R}$, $B_{2}=\{x\in \mathbb{R}:\left\Vert x\right\Vert \leq
2\}=[-2,2],$ $B_{1}=[-1,1]$ and $E(\mathbb{R},[-1,1])=\{(x-r,x+r)\subset 
\mathbb{R}:x\in \lbrack -1,1],r\geq 0\}.$
\end{example}

Let $T_{1}:[-2,2]\rightarrow 2^{E(\mathbb{R},[-1,1])}$ be defined by

$T_{1}(x)=\left\{ 
\begin{array}{c}
\lbrack -1.99995,2.00005],\text{ \ \ \ \ \ \ \ \ \ \ if \ \ \ \ \ \ \ \ \ }%
x\in \lbrack -\frac{1}{100},\frac{1}{100}]; \\ 
\lbrack -2+\frac{1}{2}x^{2},2+\frac{1}{2}x^{2}],\text{ if }x\in \lbrack -1,-%
\frac{15}{32})\cup (-\frac{15}{32},-\frac{1}{100}]\cup \\ 
\text{ \ \ \ \ \ \ \ \ \ \ \ \ \ \ \ \ \ \ \ \ \ \ \ \ \ \ \ \ \ \ \ \ \ \ \
\ \ \ }\cup (\frac{1}{100},\frac{15}{32})\cup (\frac{15}{32},1]; \\ 
\lbrack -\frac{19}{10},\frac{5}{2}],\text{ \ \ \ \ \ \ \ \ \ \ \ \ \ \ \ \ \
\ if \ \ \ \ \ \ \ \ \ \ \ \ \ \ \ }x\in \{-\frac{15}{32},\frac{15}{32}\};
\\ 
\lbrack -\frac{3}{2},\frac{5}{2}],\text{ \ \ \ \ \ \ \ \ \ \ if \ \ \ \ \ \
\ \ \ \ \ \ \ \ \ \ \ }x\in \lbrack -2,-1)\cup (1,2].%
\end{array}%
\right. $

\textit{We note that} $T_{1}$\textit{\ is almost lower semicontinuous,
condensing, with non-empty, convex and closed values.}

\textit{Let }$\Omega =[-2,2],$\textit{\ }$\mathcal{F}$ \textit{be} \textit{%
the }$\sigma -$\textit{algebra of the borelian sets of }$[-2,2]$\textit{\
and let }$T:\Omega \times \lbrack -2,2]\rightarrow E(\mathbb{R},[-1,1])$%
\textit{\ be the random operator defined by}

$T(\omega ,x)=\left\{ 
\begin{array}{c}
T_{1}(x)\text{ \ \ \ \ if \ \ \ \ \ \ }x=\omega ; \\ 
\lbrack -1.99995,\frac{5}{2}]\text{ if }x\neq \omega%
\end{array}%
\text{ }\right. $ for each $(\omega ,x)\in \Omega \times \lbrack -2,2].$

\textit{As in Example 1,} \textit{we can prove that, for each }$\omega \in
\Omega ,$ $(T(\omega ,\cdot ))^{-1}:\mathbb{R}\rightarrow 2^{[-2,2]}$\textit{%
\ is closed valued.}

\textit{Let }$\xi :[-2,2]\rightarrow \lbrack -1,1]$\textit{\ be defined by }$%
\xi (\omega )=1$\textit{\ for each }$\omega \in \lbrack -2,2].$

\textit{For each }$\omega \in \lbrack -2,2],$\textit{\ }

$T(\omega ,\xi (\omega ))=\left\{ 
\begin{array}{c}
T_{1}(1)\text{ \ \ \ \ if \ \ \ \ \ }\omega =1; \\ 
\lbrack -1.99995,\frac{5}{2}]\text{ if }\omega \neq 1.%
\end{array}%
\text{ }\right. =\left\{ 
\begin{array}{c}
\lbrack -\frac{3}{2},\frac{5}{2}]\text{\ \ \ \ if \ \ \ \ \ }\omega =1; \\ 
\lbrack -1.99995,\frac{5}{2}]\text{ if }\omega \neq 1.%
\end{array}%
\text{ }\right. $

\textit{Let }$\eta :[-2,2]\rightarrow \mathbb{R}$\textit{\ be defined by }$%
\eta (\omega )=1.00005$\textit{\ for each }$\omega \in \lbrack -2,2].$

\textit{Then, for each} $\omega \in \Omega ,$ \textit{we have }$\eta (\omega
)\in T(\omega ,\xi (\omega ))$ \textit{and}\ $0.00005=\left\Vert \eta
(\omega )-\xi (\omega )\right\Vert =d(\eta (\omega ),[-1,1])=d(\eta (\omega
),\overline{I_{[-1,1]}(\xi (\omega ))}),$ where $I_{[-1,1]}(\xi (\omega
))=I_{[-1,1]}(1)=\{1+r(y-1):y\in \lbrack -1,1],$ $r\geq 0\}=(-\infty
,1].\smallskip $

The following corollary is a random approximation result concerning the
random lower semicontinuous, condensing operators.

\begin{corollary}
\textit{Let }$(\Omega ,\mathcal{F})$\textit{\ be a measurable space, }$B_{2}$%
\textit{\ be separable in a finite dimensional Banach space }$X$\textit{\
and let }$T:\Omega \times B_{2}\rightarrow E(X,B_{1})$\textit{\ be a random
operator. Suppose that, for each }$\omega \in \Omega ,$\textit{\ }$T(\omega
,\cdot )$\textit{\ is lower semicontinuous, condensing, with non-empty
values and }$(T(\omega ,\cdot ))^{-1}:X\rightarrow 2^{B_{2}}$\textit{\ is
closed valued.}
\end{corollary}

\textit{Then, there exist a measurable mapping }$\xi :\Omega \rightarrow
B_{1}$\textit{\ and a mapping }$\eta :\Omega \rightarrow X$ \textit{such that%
} \textit{for each} $\omega \in \Omega ,$ \textit{we have}

\ \ \ \ \ \ \ \ \ \ \ \ \ \ \ \ \ \ \ \ \ \ \ \ \ \ \ \ \ \ \ \ \ \ \ \ \ \ $%
\eta (\omega )\in T(\omega ,\xi (\omega ))$

\textit{and}

\ \ \ \ \ \ \ \ \ \ \ \ \ \ \ \ \ \ \ \ \ \ \ \ \ \ \ \ \ \ $\left\Vert \eta
(\omega )-\xi (\omega )\right\Vert =d(\eta (\omega ),B_{1})=d(\eta (\omega ),%
\overline{I_{B_{1}}(\xi (\omega ))}).\medskip $

By applying Theorem 2, we obtain the following fixed point theorem
concerning the random almost lower semicontinuous and condensing operators.

\begin{theorem}
\textit{Let }$(\Omega ,\mathcal{F})$\textit{\ be a measurable space, }$B_{2}$
be \textit{separable in a finite dimensional Banach space }$X$\textit{\ and
let }$T:\Omega \times B_{2}\rightarrow E(X,B_{1})$\textit{\ be a random
operator. Suppose that, for each }$\omega \in \Omega ,$\textit{\ }$T(\omega
,\cdot )$\textit{\ is almost lower semicontinuous, condensing, non-empty
valued and }$(T(\omega ,\cdot ))^{-1}:X\rightarrow 2^{B_{2}}$\textit{\ is
closed valued.}
\end{theorem}

\textit{In addition,} \textit{for each }$\omega \in \Omega $ \textit{and} $%
x\in \partial (B_{1})\backslash T(\omega ,x),$\textit{\ }$T(\omega ,\cdot )$ 
\textit{satisfies one of the following} \textit{conditions:}

\textit{i) For each }$y\in T(\omega ,x),$ $\left\Vert y-z\right\Vert
<\left\Vert y-x\right\Vert $ for some $z\in \overline{I_{B_{1}}(x)};$

\textit{ii) For each }$y\in T(\omega ,x),$ \textit{there exists }$\lambda $ 
\textit{with }$\mathit{|}\lambda |<1$\textit{\ such that }$\lambda
x+(1-\lambda )y\in \overline{I_{B_{1}}(x)};$

\textit{iii) }$T(\omega ,x)\subseteq \overline{I_{B_{1}}(x)};$

\textit{iv) For each }$\lambda \in (0,1),$ $x\notin \lambda T(\omega ,x);$

\textit{v) For each }$y\in T(\omega ,x),$ \textit{there exists }$\gamma \in
(1,\infty )$ \textit{such that} $\left\Vert y\right\Vert ^{\gamma }-1\leq
\left\Vert y-x\right\Vert ^{\gamma };$

\textit{vi) For each }$y\in T(\omega ,x),$ \textit{there exists }$\beta \in
(0,1)$ \textit{such that} $\left\Vert y\right\Vert ^{\beta }-1\geq
\left\Vert y-x\right\Vert ^{\beta }.$

\textit{Then, }$T$\textit{\ has a random fixed point.}

\begin{proof}
According to Theorem 2, there exist a measurable mapping $\xi :\Omega
\rightarrow B_{1}$ and a mapping $\eta :\Omega \rightarrow X$ such that for
each $\omega \in \Omega ,$ we have

\ \ \ \ \ \ \ \ \ \ \ \ \ \ \ \ \ \ \ \ \ \ \ \ \ \ \ \ \ \ \ \ \ \ \ \ \ \ $%
\eta (\omega )\in T(\omega ,\xi (\omega )),$ $\xi (\omega )=r(\eta (\omega
)) $

and

\ \ \ \ \ \ \ \ \ \ \ \ \ \ \ \ \ \ \ \ \ \ \ \ \ \ \ \ \ \ $\left\Vert \eta
(\omega )-\xi (\omega )\right\Vert =d(\eta (\omega ),B_{1})=d(\eta (\omega ),%
\overline{I_{B_{1}}(\xi (\omega ))}).$

We note that $d(\eta (\omega ),\overline{I_{B_{1}}(\xi (\omega ))})>0$ for
some $\omega \in \Omega $ implies $\xi (\omega )\in \partial (B_{1})$ and $%
\left\Vert \eta (\omega )\right\Vert >1.$ Indeed, if for $\omega \in \Omega $%
, $\xi (\omega )\in $int($B_{1}),$ then, $\overline{I_{B_{1}}(\xi (\omega ))}%
=E(X,B_{1})$ and $d(\eta (\omega ),\overline{I_{B_{1}}(\xi (\omega ))}=0,$
which is a contradiction.

Further, we will prove that $T$ has a random fixed point in each of the
cases i)-vi). Then, let us assume, by way of contradiction, that there is
some $\omega \in \Omega $ such that $\xi (\omega )\notin T(\omega ,\xi
(\omega )).$

Condition i) implies that $\left\Vert \eta (\omega )-z\right\Vert
<\left\Vert \eta (\omega )-\xi (\omega )\right\Vert $ for some $z\in 
\overline{I_{B_{1}}(\xi (\omega ))},$ which contradicts the choice of $\xi .$

Condition ii) implies that there exists $\lambda $ with $|\lambda |<1$\ such
that $\lambda \xi (\omega )+(1-\lambda )\eta (\omega )\in \overline{%
I_{B_{1}}(\xi (\omega ))}.$ We obtain a contradiction, in the following way:

$\left\Vert \eta (\omega )-\xi (\omega )\right\Vert \leq \left\Vert \eta
(\omega )-(\lambda \xi (\omega )+(1-\lambda )\eta (\omega ))\right\Vert $

\ \ \ \ \ \ \ \ \ \ \ \ \ \ \ \ \ \ \ \ \ \ \ $\ =\left\Vert \lambda (\eta
(\omega )-\xi (\omega ))\right\Vert $

\ \ \ \ \ \ \ \ \ \ \ \ \ \ \ \ \ \ \ \ \ \ \ \ \ \ $=|\lambda |\left\Vert
\eta (\omega )-\xi (\omega )\right\Vert $

\ \ \ \ \ \ \ \ \ \ \ \ \ \ \ \ \ \ \ \ \ \ \ \ \ \ $<\left\Vert \eta
(\omega )-\xi (\omega )\right\Vert .$

If $T$ satisfies condition iii), then it satisfies condition ii) by letting $%
\lambda =0.$

Since, $\xi (\omega )\in \partial (B_{1}),$ condition iv) implies that for
each $\lambda \in (0,1),$ $\xi (\omega )\notin \lambda T(\omega ,\xi (\omega
))$ and then, for each $\lambda \in (0,1),$ $\xi (\omega )\neq \lambda \eta
(\omega ).$ But, we have that $\xi (\omega )=\frac{\eta (\omega )}{%
\left\Vert \eta (\omega )\right\Vert }$ and $\left\Vert \eta (\omega
)\right\Vert >1,$ which is a contradiction$.$

Condition v) implies that there exists $\gamma \in (1,\infty )$ such that $%
\left\Vert \eta (\omega )\right\Vert ^{\gamma }-1\leq \left\Vert \eta
(\omega )-\xi (\omega )\right\Vert ^{\gamma }.$ Let $\lambda _{0}=\frac{1}{%
\left\Vert \eta (\omega )\right\Vert }\in (0,1)$. Then,

$\frac{\left\Vert \eta (\omega )-\xi (\omega )\right\Vert ^{\gamma }}{%
\left\Vert \eta (\omega )\right\Vert ^{\gamma }}=(1-\lambda _{0})^{\gamma
}<1-\lambda _{0}^{\gamma }=\frac{\left\Vert \eta (\omega )\right\Vert
^{\gamma }-1}{\left\Vert \eta (\omega )\right\Vert ^{\gamma }}\leq \frac{%
\left\Vert \eta (\omega )-\xi (\omega )\right\Vert ^{\gamma }}{\left\Vert
\eta (\omega )\right\Vert ^{\gamma }}$ and therefore,

$\left\Vert \eta (\omega )-\xi (\omega )\right\Vert >\left\Vert \eta (\omega
)\right\Vert -1,$ contradicting the fact that $\left\Vert \eta (\omega )-\xi
(\omega )\right\Vert =\left\Vert \eta (\omega )\right\Vert -1,$ which is
true since $\eta (\omega )\notin B_{1}.$

In case that condition vi) is fulfilled, an argument similar to the one from
above can be done.

Consequently, in all the cases i)-vi), it remains that $\xi (\omega )\in
T(\omega ,\xi (\omega ))$ for each $\omega \in \Omega .$
\end{proof}

Now, we are establishing a random fixed point theorem for random lower
semicontinuous, condensing operators.

\begin{corollary}
\textit{Let }$(\Omega ,\mathcal{F})$\textit{\ be a measurable space, }$B_{1}$
be \textit{separable in a finite dimensional Banach space }$X$\textit{\ and
let }$T:\Omega \times B_{1}\rightarrow E(X,B_{1})$\textit{\ be a random
operator. Suppose that, for each }$\omega \in \Omega ,$\textit{\ }$T(\omega
,\cdot )$\textit{\ is lower semicontinuous, condensing, with non-empty
values and }$(T(\omega ,\cdot ))^{-1}:X\rightarrow 2^{B_{1}}$\textit{\ is
closed valued.}
\end{corollary}

\textit{In addition,} \textit{for each }$\omega \in \Omega $ and $x\in
\partial (B_{1})\backslash T(\omega ,x),$\textit{\ }$T(\omega ,\cdot )$ 
\textit{satisfies one of the following} \textit{conditions:}

\textit{i) For each }$y\in T(\omega ,x),$ $\left\Vert y-z\right\Vert
<\left\Vert y-x\right\Vert $ for some $z\in \overline{I_{B_{1}}(x)};$

\textit{ii) For each }$y\in T(\omega ,x),$ \textit{there exists }$\lambda $ 
\textit{with }$\mathit{|}\lambda |<1$\textit{\ such that }$\lambda
x+(1-\lambda )y\in \overline{I_{B_{1}}(x)};$

\textit{iii) }$T(\omega ,x)\subseteq \overline{I_{B_{1}}(x)};$

\textit{iv) For each }$\lambda \in (0,1),$ $x\notin \lambda T(\omega ,x);$

\textit{v) For each }$y\in T(\omega ,x),$ \textit{there exists }$\gamma \in
(1,\infty )$ \textit{such that} $\left\Vert y\right\Vert ^{\gamma }-1\leq
\left\Vert y-x\right\Vert ^{\gamma };$

\textit{vi) For each }$y\in T(\omega ,x),$ \textit{there exists }$\beta \in
(0,1)$ \textit{such that} $\left\Vert y\right\Vert ^{\beta }-1\geq
\left\Vert y-x\right\Vert ^{\beta }.$

\textit{Then, }$T$\textit{\ has a random fixed point.\medskip }

We introduce the following condition, which is necessary for the statement
of our next result.

\begin{definition}
condition $\mathcal{M}$:
\end{definition}

\ \ \ \ \ \ \ \ \textit{Suppose that for each }$n\in N,$\textit{\ }$\eta
_{n} $\textit{,}$\xi _{n}:\Omega \rightarrow C\subset X$\textit{\ are
measurable }

\textit{\ \ \ \ \ \ \ \ \ and for each }$\omega \in \Omega ,$\textit{\ }$%
\eta _{n}(\omega )\in T(\omega ,\xi _{n}(\omega )).$\textit{\ If for each }$%
\omega \in \Omega ,$

\textit{\ \ \ \ \ \ \ \ }$\xi _{n}(\omega )-\eta _{n}(\omega )\rightarrow 0$%
\textit{\ as }$n\rightarrow \infty ,$\textit{\ then there exists a random
fixed }

\textit{\ \ \ \ \ \ \ \ \ point }$\xi $\textit{\ for }$T.\medskip $

Now, we are obtaining a random fixed point theorem for random almost lower
semicontinuous, 1-set-contractive operators which satisfy the condition $%
\mathcal{M}$.

\begin{theorem}
\textit{Let }$(\Omega ,\mathcal{F})$\textit{\ be a measurable space, }$C$%
\textit{\ be a closed convex bounded separable subset of a finite
dimensional Banach space }$X$\textit{\ and let }$T:\Omega \times
C\rightarrow 2^{C}$\textit{\ be a random operator which satisfies condition $%
\mathcal{M}$. Let us suppose that, for each }$\omega \in \Omega ,$\textit{\ }%
$T(\omega ,\cdot )$\textit{\ is almost lower semicontinuous,
1-set-contractive with non-empty, convex and closed values and }$(T(\omega
,\cdot ))^{-1}:C\rightarrow 2^{C}$\textit{\ is closed valued. }
\end{theorem}

\ \ \ \ \ \textit{Then, }$T$\textit{\ has a random fixed point.}

\begin{proof}
We define the sequence of operators $\{T_{n}\},$ where, for each $n\in N,$ $%
T_{n}(\omega ,\cdot ):\Omega \times C\rightarrow 2^{C}$ and $T_{n}(\omega
,x)=(1-\frac{1}{n})T(\omega ,x)$ for each $\omega \in \Omega $ and $x\in C.$

We notice that for each $n\in N,$ $T_{n}(\omega ,\cdot )$\textit{\ }is
almost lower semicontinuous with non-empty, convex and closed values and $%
(T_{n}(\omega ,\cdot ))^{-1}:C\rightarrow 2^{C}$\ is closed valued.

In addition, we claim that for each $\omega \in \Omega ,$ $T_{n}(\omega
,\cdot )$ is condensing. Indeed, if we consider a bounded subset $A$ of $C$
such that $\gamma (A)>0,$ then, for each $\omega \in \Omega ,$ $\gamma
(T_{n}(\omega ,A))=(1-\frac{1}{n})\gamma (T(\omega ,A))\leq (1-\frac{1}{n}%
)\gamma (A)<\gamma (A).$

Thus, the claim is shown. All the assumptions of Theorem 2\ are fulfilled,
so, according to this result, each $T_{n}:\Omega \times C\rightarrow 2^{C}$
has a random fixed point $\xi _{n}:\Omega \rightarrow C.$ Obviously, $\xi
_{n}(\omega )\in (1-\frac{1}{n})T(\omega ,\xi _{n}(\omega ))$ for each $%
\omega \in \Omega .$

Let us consider $\eta _{n}:\Omega \rightarrow C$ such that, for each $\omega
\in \Omega ,$ $\eta _{n}(\omega )\in T(\omega ,\xi _{n}(\omega ))$ and $\xi
_{n}(\omega )=(1-\frac{1}{n})\eta _{n}(\omega ).$ Then, $\eta _{n}$ is
measurable and, since $C$ is bounded, $\xi _{n}(\omega )-\eta _{n}(\omega
)\rightarrow 0$ as $n\rightarrow \infty ,$ for each $\omega \in \Omega .$ We
can conclude that there exists $\xi :\Omega \rightarrow C$ measurable such
that $\xi (\omega )\in T(\omega ,\xi (\omega ))$ for each $\omega \in \Omega
.$
\end{proof}

A random fixed point theorem for random lower semicontinuous,
1-set-contractive operators which fulfill the condition $\mathcal{M}$ is
established now.

\begin{corollary}
\textit{Let }$(\Omega ,\mathcal{F})$\textit{\ be a measurable space, }$C$%
\textit{\ be a closed convex bounded separable subset of a finite
dimensional Banach space }$X$\textit{\ and let }$T:\Omega \times
C\rightarrow 2^{C}$\textit{\ be a random operator which satisfies condition $%
\mathcal{M}$. Let us suppose that, for each }$\omega \in \Omega ,$\textit{\ }%
$T(\omega ,\cdot )$\textit{\ is lower semicontinuous, 1-set-contractive,
with non-empty, convex and closed values and }$(T(\omega ,\cdot
))^{-1}:C\rightarrow 2^{C}$\textit{\ is closed valued. }
\end{corollary}

\ \ \ \ \ \textit{Then, }$T$\textit{\ has a random fixed point.\medskip }

A random approximation theorem for random almost lower semicontinuous,
1-set-contractive operators is established now.

\begin{theorem}
\textit{Let }$(\Omega ,\mathcal{F})$\textit{\ be a measurable space, }$B_{2}$%
\textit{\ be separable in a finite dimensional Banach space }$X$\textit{\
and let }$T:\Omega \times B_{2}\rightarrow E(X,B_{1})$\textit{\ be a random
operator which satisfies the condition $\mathcal{M}$. Let us suppose that,
for each }$\omega \in \Omega ,$\textit{\ }$T(\omega ,\cdot )$\textit{\ is
almost lower semicontinuous, 1-set-contractive, with non-empty values, and }$%
(T(\omega ,\cdot ))^{-1}:X\rightarrow 2^{B_{2}}$\textit{\ is closed valued.}
\end{theorem}

\textit{Then, there exist a measurable mapping }$\xi :\Omega \rightarrow
B_{1}$\textit{\ and a mapping }$\eta :\Omega \rightarrow X$ \textit{such that%
} \textit{for each} $\omega \in \Omega ,$ \textit{we have}

\ \ \ \ \ \ \ \ \ \ \ \ \ \ \ \ \ \ \ \ \ \ \ \ \ \ \ \ \ \ \ \ \ \ \ \ \ \ $%
\eta (\omega )\in T(\omega ,\xi (\omega ))$

\textit{and}

\ \ \ \ \ \ \ \ \ \ \ \ \ \ \ \ \ \ \ \ \ \ \ \ \ \ \ \ \ \ $\left\Vert \eta
(\omega )-\xi (\omega )\right\Vert =d(\eta (\omega ),B_{1})=d(\eta (\omega ),%
\overline{I_{B_{1}}(\xi (\omega ))}).$

\begin{proof}
Let us define $r:X\rightarrow B_{1}$ by $r(x)=\left\{ 
\begin{array}{c}
x\text{ if }x\in B_{1}; \\ 
\frac{x}{\left\Vert x\right\Vert }\text{ if }x\notin B_{1}.%
\end{array}%
\right. $

Then, $r$ is continuous and $r(A)\subseteq \overline{\text{co}}(A\cup \{0\})$
for each bounded subset $A$ of $X.$ In addition, for each $\omega \in \Omega
,$ $T(\omega ,\cdot )$ is 1-set-contractive$,$ and we conclude that $%
G(\omega ,\cdot )=r\circ T(\omega ,\cdot ):B_{2}\rightarrow 2^{B_{1}}$ is
1-set-contractive. According to the hypotheses, $T(\omega ,\cdot )$\ is
almost lower semicontinuous with non-empty convex closed values and $%
(T(\omega ,\cdot ))^{-1}:X\rightarrow 2^{B_{1}}$\ is closed valued. Hence,
it is easy to check that $G(\omega ,\cdot )$\ is almost lower semicontinuous
with non-empty convex closed values, satisfies the condition \textit{$%
\mathcal{M}$} and $(G(\omega ,\cdot ))^{-1}:B_{1}\rightarrow 2^{B_{2}}$\ is
closed valued. Then, $G$ fulfills all the conditions of Theorem 4. By
applying this theorem, we obtain that there exists $\xi :\Omega \rightarrow
B_{1}$ such that for each $\omega \in \Omega ,$ we have $\xi (\omega )\in
G(\omega ,\xi (\omega ))$. Further, the proof follows the same line as the
proof of Theorem 2.
\end{proof}

The result below is a random approximation theorem for random lower
semicontinuous, 1-set-contractive operators which satisfy the condition $%
\mathcal{M}$.

\begin{corollary}
\textit{Let }$(\Omega ,\mathcal{F})$\textit{\ be a measurable space, }$B_{2}$%
\textit{\ be separable in a finite dimensional Banach space }$X$\textit{\
and let }$T:\Omega \times B_{2}\rightarrow E(X,B_{1})$\textit{\ be a random
operator which satisfies the condition $\mathcal{M}$. Suppose that, for each 
}$\omega \in \Omega ,$\textit{\ }$T(\omega ,\cdot )$\textit{\ is lower
semicontinuous, 1-set-contractive, with non-empty values and }$(T(\omega
,\cdot ))^{-1}:X\rightarrow 2^{B_{2}}$\textit{\ is closed valued.}
\end{corollary}

\textit{Then, there exist a measurable mapping }$\xi :\Omega \rightarrow
B_{1}$\textit{\ and a mapping }$\eta :\Omega \rightarrow X$ \textit{such that%
} \textit{for each} $\omega \in \Omega ,$ \textit{we have}

\ \ \ \ \ \ \ \ \ \ \ \ \ \ \ \ \ \ \ \ \ \ \ \ \ \ \ \ \ \ \ \ \ \ \ \ \ \ $%
\eta (\omega )\in T(\omega ,\xi (\omega ))$

\textit{and}

\ \ \ \ \ \ \ \ \ \ \ \ \ \ \ \ \ \ \ \ \ \ \ \ \ \ \ \ \ \ $\left\Vert \eta
(\omega )-\xi (\omega )\right\Vert =d(\eta (\omega ),B_{1})=d(\eta (\omega ),%
\overline{I_{B_{1}}(\xi (\omega ))}).\medskip $

By using the above result, we state the existence of the random fixed points
for random almost lower semicontinuous, 1-set-contractive operators under
the following assumptions.

\begin{theorem}
\textit{Let }$(\Omega ,\mathcal{F})$\textit{\ be a measurable space, }$B_{2}$
be \textit{separable in a finite dimensional Banach space }$X$\textit{\ and
let }$T:\Omega \times B_{2}\rightarrow E(X,B_{1})$\textit{\ be a random
operator which satisfies condition $\mathcal{M}$. Let us suppose that, for
each }$\omega \in \Omega ,$\textit{\ }$T(\omega ,\cdot )$\textit{\ is almost
lower semicontinuous, 1-set-contractive, with non-empty values and }$%
(T(\omega ,\cdot ))^{-1}:X\rightarrow 2^{B_{2}}$\textit{\ is closed valued.}
\end{theorem}

\textit{In addition, let us suppose that for each }$\omega \in \Omega $ and $%
x\in \partial (B_{1})\backslash T(\omega ,x),$\textit{\ }$T(\omega ,\cdot )$ 
\textit{satisfies one of the following} \textit{conditions:}

\textit{i) For each }$y\in T(\omega ,x),$ $\left\Vert y-z\right\Vert
<\left\Vert y-x\right\Vert $ for some $z\in \overline{I_{B_{1}}(x)};$

\textit{ii) For each }$y\in T(\omega ,x),$ \textit{there exists }$\lambda $ 
\textit{with }$\mathit{|}\lambda |<1$\textit{\ such that }$\lambda
x+(1-\lambda )y\in \overline{I_{B_{1}}(x)};$

\textit{iii) }$T(\omega ,x)\subseteq \overline{I_{B_{1}}(x)};$

\textit{iv) For each }$\lambda \in (0,1),$ $x\notin \lambda T(\omega ,x);$

\textit{v) For each }$y\in T(\omega ,x),$ \textit{there exists }$\gamma \in
(1,\infty )$ \textit{such that} $\left\Vert y\right\Vert ^{\gamma }-1\leq
\left\Vert y-x\right\Vert ^{\gamma };$

\textit{vi) For each }$y\in T(\omega ,x),$ \textit{there exists }$\beta \in
(0,1)$ \textit{such that} $\left\Vert y\right\Vert ^{\beta }-1\geq
\left\Vert y-x\right\Vert ^{\beta }.$

\textit{Then, }$T$\textit{\ has a random fixed point.}

\begin{proof}
This result is an application of Theorem 5. The proof follows the same line
as the proof of Theorem 3.\bigskip
\end{proof}

By using the above result, we finally state the existence of the random
fixed points for the random lower semicontinuous, 1-set-contractive
operators under the following assumptions.

\begin{corollary}
\textit{Let }$(\Omega ,\mathcal{F})$\textit{\ be a measurable space, }$B_{2}$
be \textit{separable in a finite dimensional Banach space }$X$\textit{\ and
let }$T:\Omega \times B_{2}\rightarrow E(X,B_{1})$\textit{\ be a random
operator which satisfies the condition $\mathcal{M}$. Let us suppose that,
for each }$\omega \in \Omega ,$\textit{\ }$T(\omega ,\cdot )$\textit{\ is
lower semicontinuous, 1-set-contractive, with non-empty values and }$%
(T(\omega ,\cdot ))^{-1}:X\rightarrow 2^{B_{2}}$\textit{\ is closed valued.}
\end{corollary}

\textit{If,} \textit{for each }$\omega \in \Omega $ and $x\in \partial
(B_{1})\backslash T(\omega ,x),$\textit{\ }$T(\omega ,\cdot )$ \textit{%
satisfies one of the following} \textit{conditions:}

\textit{i) For each }$y\in T(\omega ,x),$ $\left\Vert y-z\right\Vert
<\left\Vert y-x\right\Vert $ for some $z\in \overline{I_{B_{1}}(x)};$

\textit{ii) For each }$y\in T(\omega ,x),$ \textit{there exists }$\lambda $ 
\textit{with }$\mathit{|}\lambda |<1$\textit{\ such that }$\lambda
x+(1-\lambda )y\in \overline{I_{B_{1}}(x)};$

\textit{iii) }$T(\omega ,x)\subseteq \overline{I_{B_{1}}(x)};$

\textit{iv) For each }$\lambda \in (0,1),$ $x\notin \lambda T(\omega ,x);$

\textit{v) For each }$y\in T(\omega ,x),$ \textit{there exists }$\gamma \in
(1,\infty )$ \textit{such that} $\left\Vert y\right\Vert ^{\gamma }-1\leq
\left\Vert y-x\right\Vert ^{\gamma };$

\textit{vi) For each }$y\in T(\omega ,x),$ \textit{there exists }$\beta \in
(0,1)$ \textit{such that} $\left\Vert y\right\Vert ^{\beta }-1\geq
\left\Vert y-x\right\Vert ^{\beta },$

\textit{then, }$T$\textit{\ has a random fixed point.}

\section{CONCLUDING REMARKS}

We have proved the existence of random approximation and fixed points for
almost lower semicontinuous and lower semicontinuous operators defined on
finite dimensional Banach spaces. Our study extends on some results which
exist in literature. It is an interesting problem which deserves further
research as to establish new similar theorems for other types of operators.

\end{document}